\documentstyle[12pt,twoside]{article}

\date{}

\begin{document}


\centerline{}

\centerline{}

\centerline {\large{\bf A lattice-ordered skew-field is totally
ordered}}

\centerline{\large{\bf if squares are positive}}
\bigskip

\centerline{\bf Yang, Yichuan}

\centerline{}

\centerline{Department of Mathematics, Beihang University,}
\centerline{100083, Beijing, P. R. China} 
\centerline{E-mail: yichuanyang@hotmail.com}

\newtheorem{Theorem}{\quad Theorem}[section]

\newtheorem{Definition}[Theorem]{\quad Definition}
\newtheorem{Question}[Theorem]{\quad Question}

\newtheorem{Corollary}[Theorem]{\quad Corollary}

\newtheorem{Lemma}[Theorem]{\quad Lemma}

\newtheorem{Remark}[Theorem]{\quad Remark}

\newtheorem{Example}[Theorem]{\quad Example}

\begin{abstract}
We show that a lattice-ordered field (not necessarily commutative)
is totally ordered if and only if each square is positive,
answering a generalized question of Conrad and Dauns \cite{CD}
 in the affirmative. As a
consequence, any lattice-ordered skew-field in \cite{Bru}
 is totally ordered. Furthermore,
we note that each lattice order determined by a {\it pre-positive
cone} $P$ on a skew-filed $F$ is linearly ordered since
$F^2\subseteq P$ (see \cite{P}).
\end{abstract}



\section{Introduction}
In \cite{AS}, Artin and Schreier observed that a totally ordered
commutative field cannot have negative squares, and Johnson
\cite{Jo} and Fuchs  \cite{Fu2} extended this result to totally
ordered domains with unit element.  In \cite{S}, Schwartz showed
that an Archimedean lattice-ordered (commutative) field that has
$1>0$ and that is algebraic over its maximal totally ordered
subfield cannot have negative squares, and in \cite{DS}, DeMarr
and Steger showed that in a partially ordered  finite dimensional
real linear algebra no square can be the negative of a strong
unit. Furthermore,    in \cite{Y1}, we guarantee the existence of
directed commutative fields with negative squares. In 1969, Conrad
and Dauns \cite{CD} raised the following problem (this is Question
(b) of their list in \cite{CD}, p.  397]).

  \begin{Question}\label{cdquestionbhao}
  Is it true that a lattice-ordered (commutative) field  $F$ in which each square is positive must be  totally ordered?
\end{Question}
 In fact, an affirmative answer for commutative case can be found
in Bourbaki \cite{Bo2} (Chapitre VI, p.  43) in 1952,  in Birkhoff
and Pierce
\cite{BP} (p.  59) in 1956, and in Fuchs \cite{Fu} in 1963. 
It should also be noted
 Redfield  answered  Conrad and Dauns' question independently in 1975 (see \cite{Red}, p.  124)
. However, none of these previous authors stated the result for
skew-fields.

  In \cite{Y},  we give a positive answer to  the generalized problem
of Conrad and Dauns for skew-fields
.  As a   consequence,  the
   ``lattice-ordered skew-fields'' in Brumfiel \cite{Bru} are  in fact   totally ordered since
   each square of an element is positive according to the definition (see \cite{Bru},
   p.  32). Furthermore, we note that every lattice order determined by a
   {\it pre-positive  cone} $P$ on a skew-field $F$ is  linearly ordered since
   $F^2\subseteq P$ (see Prestel \cite{P}). In this note we archived the promised note and give
   another interesting corollary (Corollary \ref{cd}) which  was
   deleted  in \cite{Y} from the original manuscript of \cite{Y} for the wide spectrum of the general $Monthly$ readership,
   since it is likely to be understandable and of interest to only specialists in ordered algebraic structures.

\section{Main results}

\begin{Lemma}\label{lemcd} {\it If $R$ is  an $l$-ring in which each  strictly positive element of $R$ is invertible and in which $a^2\geq 0$ for every $a$ in $R$, then $a^+|a|^{-1}a^-=0$ holds for all   $a$ in $R\setminus\{0\}$.}
\end{Lemma}

{\it Proof.}  Let $a\in R$, $a\not=0$, and define $b=
a^+|a|^{-1}a^+.$ By hypothesis, $(|a|^{-1})^2|a|=|a|^{-1}>0$ and
so $b\geq 0.$ From $a= a^++a^-,\ |a|=a^+-a^-$ we obtain that
 $b-a$ 
  $=  (-a^-)|a|^{-1}(-a^-)$
    $\geq 0.$
  That is, $b\geq a$, and hence $b\geq a^+.$ On the other hand, we have  
    $a^+-b$
                         $ = a^+|a|^{-1}a^-$
                            $\geq 0$.
                               This means that $a^+\geq b,$ and thus it shows that $b=a^+.$
                               Consequently, the last inequality implies 
                                $a^+|a|^{-1}a^-=0.$\hfill $\Box $ 

                               \medskip

The following theorem characterizes lattice-ordered rings which
are totally ordered skew-fields and provides a solution to
generalized Conrad and Dauns' problem.

 \begin{Theorem}\label{thmcd} If $R$ is  an $l$-ring, then the  following statements are equivalent:

 (a) each strictly positive element in $R$ has a multiplicative inverse and $a^2\geq 0$ for every  $a$   in $R$;

 (b) $R$ is a totally ordered skew-field.
  \end{Theorem}

  {\it Proof.} It suffices to prove (a)$\Rightarrow $(b): Lemma \ref{lemcd}
  above shows that  $a^+|a|^{-1}a^-=0$ holds for each nonzero element $a$ in
  $R$. Hence  $a^+=0$ or $a^-=0$ ( i.~e.   $R$ is totally ordered). Hence $R$ is a skew-field.
  \hfill $\Box$
  \smallskip

We note that it is an easy exercise for readers to verify that any
known proof for commutative case can be modified or used to prove
Theorem \ref{thmcd}. For a lattice-ordered skew-field, we extract
from Theorem \ref{thmcd} the following information:

\smallskip

\begin{Corollary}\label{cd} {\it Let $R$ be a lattice-ordered division ring. Then the following statements are equivalent:

 $(\alpha)$ $a^2\geq 0$ for all $a\in R$.

$(\beta)$ $R$ is order division-closed (that is, for all $a, \
b\in R$, $ab>0$ and one of $a,\ b$ is $>0$, then so is the other).

$(\gamma)$ $R$ is totally ordered.

$(\delta)$ if $a\in R$, then there exists a natuaral number $n_a$
such that $a^{n_a}\geq 0$.

$(\epsilon)$ $R$ is an f-ring (that is, $a\wedge b=0$ and $c\geq
0$ implies $ca\wedge b=ac\wedge b=0$).

$(\zeta)$ $R$ is an almost f-ring (that is, $a\wedge b=0$ implies
$ab=0$).

$(\lambda)$ the additive group of $R$ is the  group of
divisibility of a valuation domain.

  $(\eta)$  each strictly positive element in $R$ has a strictly positive multiplicative inverse.
 }
 \end{Corollary}

 {\it Proof.} ($\gamma$)$\Leftrightarrow (\alpha)$ is clear.
 ($\gamma$)$\Leftrightarrow (\lambda)$ follows from the well known theorem of Krull.
  For the rest of the proof, it suffices to prove $\zeta)\Leftrightarrow \epsilon)$,
  but this is a directed corollary of Lemma 1 (that is, in any $f$-ring, $a\wedge b=0$ implies $ab=0$) on p.  404
   of \cite{Bi} and the definitions of  almost $f$-rings and division rings. Finally, the equivalence of
   $(\eta)$ and $(\gamma)$ is plain by the proof of the Lemma \ref{lemcd}. \hfill $\Box$
 \smallskip

 Especially, for a lattice-ordered (commutative) field, we recover
 \cite{BP, Bo2, Fu, Red}:

\begin{Corollary} {\bf (Bourbaki,  Birkhoff and Pierce,  Fuchs, Redfield)} {\it If $R$ is a lattice-ordered (commutative) field in which every square is positive, then $R$ is totally ordered.}
\end{Corollary}


\end{document}